\documentclass[a4paper,12pt]{article}
\usepackage{amssymb,latexsym,amsmath,amsthm,amscd,graphpap}
\usepackage{amsmath}
\usepackage{amsfonts}
\usepackage{amssymb}
\usepackage[dvips]{epsfig}

\newtheorem{theorem}{Theorem}[section]
\newtheorem{prop}[theorem]{Proposition}

\newtheorem{lemma}[theorem]{Lemma}

\newcommand{\Rm}{\mathbb{R}}

\newcommand{\Lm}{\mathbb{L}}

\newcommand{\Lc}{\mathcal{L}}
\newcommand{\Pc}{\mathcal{P}}
\newcommand{\Nm}{\mathbb{N}}

\newcommand{\Com}{\mathbb{C}}

\newcommand{\Gm}{\mathfrak{G}}

\newcommand{\Cm}{\mathcal{C}}
\newcommand{\Nc}{\mathcal{N}}

\newcommand{\Om}{\mathcal{O}}
\newcommand{\Em}{\mathcal{E}}

\newcommand{\Um}{\mathcal{U}}
\newcommand{\Vm}{\mathcal{V}}

\newcommand{\Zr}{\mathbb{Z}}

\newcommand{\no}{n$^{\text{o}}$}

\def\HB {\hfill\break}

\theoremstyle{definition}
\newtheorem{definition}[theorem]{Definition}

\theoremstyle{remark}
\newtheorem{remark}[theorem]{Remark}
\newtheorem{rems}[theorem]{Remarks}

\numberwithin{equation}{section}
\def\HB {\hfill\break}

\newcommand{\Bc}{\mathcal{B}}

\begin{document}

     \title{\bf Generic hyperbolicity of equilibria and periodic orbits of
the parabolic equation on the circle}

\author{Romain JOLY and Genevi\`eve RAUGEL}
\date{}

\maketitle

\begin{abstract}
In this paper, we show that, for scalar reaction-diffusion
equations on the circle $S^1$, the property of hyperbolicity of all
equilibria and
periodic orbits is generic with respect to the non-linearity .
In other words, we
prove that in an appropriate functional space of nonlinear terms
in the equation, the set of functions, for which all equilibria and
periodic orbits are hyperbolic, is a countable
intersection of open dense sets. The main tools in the proof are the
property of
the lap number and the Sard-Smale theorem.
\end{abstract}

\section{Introduction}
In the local study of the qualitative properties of perturbed
dynamical systems, the concepts of non degeneracy and
hyperbolicity of equilibria and periodic orbits play a crucial role.
Indeed, if one slightly perturbs a continuous dynamical system whose
equilibria and periodic orbits are non-degenerate, then at least, the
perturbed system still admits equilibria and periodic orbits nearby.
If these elements are in addition hyperbolic, then no local
bifurcation phenomena occur, which means that the dynamics in the
neighborhood of hyperbolic equilibria and periodic orbits is stable
under small perturbations.  This stability property has practical
consequences.  For example, it implies that a numerical simulation of
the dynamics near a hyperbolic equilibrium or periodic orbit is
qualitatively correct.  These considerations show that it is important
to know if, in a given class of dynamical systems or evolutionary
equations, the equilibria and periodic orbits are all hyperbolic or
if, at least, this property is generically satisfied.\HB
Generic
hyperbolicity of equilibrium points and periodic orbits is well-known
in the case of diffeomorphisms or vector fields on a compact
finite-dimensional manifold (for instance, see Chapter 3 of
\cite{PalisMelo} and the references therein).  The genericity of
hyperbolicity of equilibria and periodic orbits was generalized in
1977 to the case of functional differential equations defined in
$\Rm^n$, $n \geq 1$, by Mallet-Paret (\cite{Mallet}).  In the frame of
evolutionary partial differential equations, generic hyperbolicity of
equilibrium points is also a classical result (see for example,
\cite{Babin-V-livre}, \cite{Bruno-Chow}, \cite{Bruno-Pola},
\cite{SW} for parabolic equations).  Notice that this is sufficient
for gradient dynamical systems, where no periodic orbits can occur.
On the contrary, for autonomous evolutionary partial differential equations,
which are not of gradient type, the hyperbolicity of periodic
orbits seems to be an open question.
\vskip 1mm

In this paper, we address the question of genericity of hyperbolicity of
periodic orbits by considering one of the simplest non-gradient
partial differential equations, namely the following scalar
reaction-diffusion equation on
the one-dimensional torus (or unit circle) $S^1= \Rm / 2 \pi
\Zr$, given by
\begin{equation}\label{eq}
\begin{split}
u_t(x,t) &=u_{xx}(x,t)+f(x,u(x,t),u_x(x,t))~, \quad (x,t)\in S^1\times\Rm_+^*~,
\cr
u(x,0) & = u_0(x)~, \quad x\in S^1~,
\end{split}
\end{equation}
where $f$ belongs to the space $C^2(S^1\times \Rm\times\Rm,\Rm)$ and
$u_0$ is given in  the Sobolev space $H^s(S^1)$, with $s\in(3/2,2)$
(so that $H^s(S^1)$ is continuously embedded into $C^{1 + \alpha}(S^1)$
for $\alpha=s -3/2$). Our purpose is to prove that the equilibria and
the periodic orbits of \eqref{eq} are hyperbolic generically with
respect to $f$. We will also consider the more constricting case
where the non-linearity is independent of $x$
\begin{equation}\label{eq-inde-x}
u_t(x,t)=u_{xx}(x,t)+f(u(x,t),u_x(x,t))~, \quad(x,t)\in S^1\times\Rm_+^*~.
\end{equation}

Equations \eqref{eq} and \eqref{eq-inde-x} are good problems to begin
with. Indeed, on
one hand these equations may have periodic orbits (see Theorem 2 of
\cite{Sandstede-Fiedler}) and on the other hand, they are the
simplest non-gradient PDE's to study since they admit very special
properties as explained below.
\vskip 1mm

In the study of global stability, generic hyperbolicity of closed
orbits plays also an important role. Showing genericity of closed
orbits is the first step in the proof the genericity of Morse-Smale property.
Morse-Smale dynamical systems are systems whose non-wandering set
consists only in a finite number of hyperbolic equilibria and
hyperbolic periodic orbits and for which the intersections of the
stable and unstable manifolds of equilibria and periodic orbits are
all transversal. On finite-dimensional compact manifolds, the
dynamics of Morse-Smale dynamical systems are globally stable under
perturbations, that is, any small regular enough perturbation of a
Morse-Smale system is still Morse-Smale and the flows are
topologically equivalent (for more details, see \cite{PalisMelo}).
     In the class of gradient dynamical systems on finite-dimensional compact
     manifolds, the Morse-Smale property is known to be generic. For
     the non-gradient dynamical systems, the Morse-Smale property is
     generic only if the manifold
is one- or two-dimensional. In the infinite-dimensional case, the structural
stability of Morse-Smale systems remains true under compactness
assumptions (for instance, under the assumption that the dynamical
systems admit compact global attractors). The genericity of
Morse-Smale property has been obtained for some classes of gradient
partial differential equations (\cite{Ang86, Henry85, Bruno-Pola,
Bruno-Raugel, RJ0}).
For Morse-Smale non-gradient evolutionary partial differential
equations, Equation \eqref{eq} is the right candidate since its asymptotic
behaviour in time is analogous to the one of
a system of ordinary differential equations in $\Rm^2$. Indeed, in
\cite{FMP}, Fiedler and Mallet-Paret have proved that Equation
\eqref{eq} satisfies a Poincar\'{e}-Bendixson type property.  For
non-linearities independent of $x$, automatic transversality
properties have been proved in \cite{FRW}, by showing that, for
\eqref{eq-inde-x}, the Morse-Smale property is equivalent to the
hyperbolicity of all closed orbits.  The results of the present paper
thus imply the genericity of Morse-Smale property for Equation
\eqref{eq-inde-x}.  In the more general case of Equation \eqref{eq},
R. Czaja and C. Rocha have shown in \cite{CR} that, if the periodic
orbits of \eqref{eq} are all hyperbolic, then their stable and
unstable manifolds intersect transversally.  Here we will show that
this hyperbolicity is generic in the non-linearity $f(x, u, u_x)$.
This shows, together with the results of \cite{CR}, that the
transversality of the stable and unstable manifolds of the periodic orbits
of \eqref{eq} is generic with respect to $f$.  Thus, this paper is
also a further step towards the proof of genericity of the Morse-Smale
property for \eqref{eq}.  We shall conclude this proof in the
forthcoming paper \cite{JR}.
\vskip 1mm

\indent When the non-linearity $f$ does not explicitly
depend on $x$ and is dissipative, the dynamics of
Equation \eqref{eq-inde-x} can be described with more precision. In
\cite{Angenent-Fiedler}, Angenent and Fiedler obtained the following
     characterization of the closed orbits of \eqref{eq-inde-x}.

\begin{prop}\label{th-Ang-Fied}
Any closed orbit $u(x,t)$ of \eqref{eq-inde-x} is:\HB
- either a homogeneous equilibrium point, that is,
$u(x,t)\equiv e\in\Rm$ where $f(e,0)=0$,\\
- or a rotating wave, that is, $u(x,t)=v(x-ct)$ with $c\neq 0$, \HB
- or a frozen wave, that is, $u$ is a non-homogeneous solution of
$u_{xx}+f(u,u_x)=0$.
\end{prop}

Moreover, Angenent and Fiedler have proved that the $\omega$-limit
set of any element $u_0 \in H^s(S^1)$ contains a rotating wave or a
steady state.
Under the hypothesis that the rotating waves and
equilibria are all hyperbolic, they showed that the $\omega$-limit
set consists exactly in one rotating wave or in one equilibrium
point and also described the connecting orbits.
In \cite{MN}, Matano and Nakamura proved that there exist
neither homoclinic orbits nor heteroclinic cycles.
In \cite{FRW}, Fiedler, Rocha and Wolfrum, using the notion of
``adjacency", have characterized
connecting orbits under the hypothesis that the rotating waves and
     equilibria are hyperbolic. They also gave an equivalent
     characterization of the property of hyperbolicity of rotating
waves. Thus, the generic hyperbolicity of the closed orbits of
\eqref{eq-inde-x} is primordial to conclude that the dynamics
described in these papers include every possible dynamics, except
maybe some exceptional ones.\\

Finally, let us recall that, if the
periodic boundary conditions are replaced by homogeneous
Dirichlet or Neumann boundary ones, then the dynamical system
$S(t)$ associated with \eqref{eq} is of gradient type
and  the $\omega$-limit set of any point is exactly one
equilibrium point (see \cite{Ze}). Moreover, in these cases, as shown
by D. Henry \cite{Henry85} and later by S. Angenent \cite{Ang86}, the
stable and unstable manifolds of the equilibria intersect
transversally, provided they are all hyperbolic. Thus, the Morse-Smale
property is generic in the class of one-dimensional scalar parabolic
equations on the segment with Dirichlet or Neumann boundary
conditions. In a series of
papers, Fiedler and Rocha have described the different equivalence
classes of
attractors and, in particular, the heteroclinic orbits connecting the
equilibria (see \cite{FR96}, \cite{FR00}, for example).

\vskip 2mm

Now we are going to precisely state the genericity results proved in
this
paper. We set $\Gm=\Cm^2(S^1\times \Rm\times\Rm,\Rm)$ and we endow
$\Gm$ with the Whitney
topology, that is, the topology generated by the neighborhoods
\begin{equation}\label{2-topo-whitney}
\{g\in\Gm~/~|D^i f(x,u,v)-D^i g(x,u,v)|\leq\delta(u,v),~\forall
i\in\{0,1,2\},~
\forall (x,u,v)\in S^1\times\Rm^2\}~,
\end{equation}
where $f$ is any function in $\Gm$ and $\delta$ is any positive
continuous function
(see \cite{GG}). \HB
It is well known that $\Gm$ is a Baire space, which means that any
generic set, that is any countable intersection of open and dense
sets, is dense in $\Gm$ (see \cite{GG} for instance).

In the general case, when $f$ depends also on the variable $x$, we
prove the following genericity result (for the definition of
hyperbolicity we refer the reader to Section \ref{Sec-spectre} below).

\begin{theorem} \label{Peixoto}
There exists a generic subset $\Om$ of $\Gm$ such that, for all $f
\in  \Om$, all the equilibrium points and all the periodic solutions
of \eqref{eq} are hyperbolic.
\end{theorem}

\begin{rems}\HB
1) The above theorem is also valid if, in the equation \eqref{eq}, we
replace the Laplacian $-\partial_{xx}$ by the self-adjoint operator
$- a(x) \partial_{xx}$ or $-a(x)^{-1} (\partial_x( a(x) \partial_{x}
\cdot))$, where $a(x) >0$ is a $C^1$-function on $S^1$. \HB
2) In Theorems \ref{Peixoto} and \ref{th-inde-x}, we assume for
sake of simplicity that $f$ is of class $C^2$ in all the variables.
Actually, we can begin with a nonlinearity $f$, which  is less
regular, and immediately replace it by a small perturbation $f_0$,
which is in $C^2(S^1\times \Rm\times\Rm,\Rm)$.

\end{rems}

\vskip2mm

As we already mentioned, it is also interesting to consider a
non-linearity $f$, which is independent of $x$. Even if the choice of
such non-linearity is more
constrained, this case is not more difficult than the general case
due to the particular properties coming from the $S^1$-equivariance
of \eqref{eq-inde-x}. In Section \ref{Sec-inde}, we prove the
following genericity result.

\begin{theorem}\label{th-inde-x}
Let $\Gm^i$ be the subspace of $\Gm$ consisting of functions
independent of $x$. There exists a generic subset $\Om^i$ of $\Gm^i$
such that, for all
$f\in\Om^i$, \eqref{eq-inde-x} has no frozen waves and any
homogeneous equilibrium or rotating wave of \eqref{eq-inde-x} is
hyperbolic.
\end{theorem}

\begin{remark}
     In \cite{FRW}, the generic hyperbolicity of homogeneous equilibria and
rotating waves is asserted, based on a personal communication of
P. Brunovsk\'{y}.
\end{remark}
\vskip 4mm

There are two  main ingredients in the proof of Theorems
\ref{Peixoto} and
\ref{th-inde-x}.
The first one, as expected in genericity properties, is
the Sard-Smale Theorem that we recall in the appendix. The second main
ingredient is the Sturm property (also called zero number or lap
number property).
This property originated in the paper of Sturm \cite{Sturm} (see also
\cite{Matano82} and \cite{Angenent} for more recent results, mainly in
the case where $f$ is not analytic). For any $\varphi \in C^1(S^1)$,
we define the zero number $z(\varphi)$ as the (even) number of strict
sign changes of $\varphi$. Let $I$ be an open subinterval of $\Rm$
and let $v(x,t)$ be a solution
of the linear parabolic equation
\begin{equation}
\label{parabolic}
v_t = v_{xx} + b(x,t) v + d(x,t) v_x~,
\end{equation}
where $b$, $b_x$, $b_t$ and $d$ are bounded functions on any compact
subset of $S^1 \times I$. Then $z(v(\cdot,t))$
is finite, for any $t >0$, and nonincreasing with $t$. Moreover, the
zero number $z(v(\cdot,t))$ drops strictly at $t=t_0$,
if and only if there exists $x_0 \in S^1$ such that
\begin{equation}
\label{double}
v(x_0,t_0)=0~, \quad \partial_x v(x_0,t_0)=0~.
\end{equation}
In Section 2.3, we will see that these remarkable properties of the
zero number imply also special spectral properties for the linearized
equation associated with \eqref{eq}. In the proof of Theorem
\ref{Peixoto}, it will also play an important role. \HB
Finally, we point out that the zero number property has been widely
used in
all the above mentioned papers dealing with the description of the
dynamics of a scalar one-dimensional equation.
\vskip 3mm
The paper is organized as follows. In Section 2, we give the main
spectral properties of the linearized equations around the equilibria
and the periodic solutions of \eqref{eq}. Section 3 contains
the proof of the genericity of the hyperbolicity of the equilibria of
\eqref{eq}, while Section 4 is devoted to the genericity of the
hyperbolicity of the periodic orbits of \eqref{eq}.  In Section 5, we
consider the special case of equation \eqref{eq-inde-x} and prove
Theorem \ref{th-inde-x}.
Finally, in the appendix, we recall the Sard-Smale theorem, that we
apply here.

\vskip 3mm

{\noindent \bf Acknowledgements:} The authors wish to thank P.
Brunovsk\'y, R. Czaja and C. Rocha for fruitful discussions with them.


\section{Spectrum of the linearized operators}\label{Sec-spectre}

The notions of non-degeneracy and hyperbolicity of equilibria and
periodic orbits of
an evolutionary equation are related to the spectral properties of
the linearized equations around these elements. In this section, we
will recall these notions and describe several properties of the
linearized equations associated with \eqref{eq}.

Let $p$ be either an equilibrium point $e$ or a
periodic orbit $p_0(x,t)$ of \eqref{eq}.
We consider the linearized equation
\begin{equation}
\label{linearized}
\begin{split}
\varphi_t &=\varphi_{xx}+ D_{u}f(x,p,p_{x})\varphi
+ D_{u_x}f(x,p,p_{x})\varphi_x~, \cr
\varphi (x,0) &= \varphi_0(x)~.
\end{split}
\end{equation}
Let $s\in (3/2,2)$; we introduce the operator $U(t,0): H^s(S^1)
\longrightarrow H^s(S^1)$,
     defined by $U(t,0)\varphi_0= \varphi(t)$ where
$\varphi(t)$ is the solution of the linearized equation
\eqref{linearized}.
Sometimes, if
needed, we will also denote the operator $U(t,0)$ by $U_{p,f}(t,0)$.

\subsection{The case of equilibria}

Let $e \in H^s(S^1)$ be an equilibrium point of \eqref{eq}, that is
a solution of $e_{xx}+f(x,e,e_x)=0$.
Due to the classical elliptic estimates and the regularity of $f$,
the equilibrium $e$ belongs to
$H^{s+2}(S^1)$. In the case where $p=e$ is an equilibrium point, the
linearized operator $U(t,0)$ is an analytic semigroup generated by
the linear  operator
\begin{equation}
\label{MapLe}
L_e=\partial_{xx}^2+f'_u(x,e,e_x)+f'_{u_x}(x,e,e_x)\partial_x :H^{2}(S^1)
\rightarrow L^2(S^1)~,
\end{equation}
that is, $U(t,0) = e^{L_e t}$.\\
The general definitions of simplicity and hyperbolicity are as
follows.
\begin{definition}$~$\\
1) We say that an equilibrium point $e$ is simple if $0$ belongs to
the resolvent set of the operator
$L_e$. \HB
2) We say that $e$ is a hyperbolic equilibrium point if the
intersection of the spectrum $\sigma(U(t,0))$ of $U(t,0)$ with the
unit circle $S^1$
in the complex plane is the empty set.
\end{definition}

Since $L_e$ is a Fredholm operator and also a sectorial operator,
simple characterizations of the notions of simplicity and
hyperbolicity can be given (for the properties of Fredholm operators,
see for example \cite{Bonic}).

\begin{lemma}\label{lemme-Le-Fredholm}
Let $e$ be any equilibrium point of \eqref{eq}.\\
1) The  operator $L_e:H^{2}(S^1)\longrightarrow L^2(S^1)$ is a
Fredholm operator of index $0$. As a consequence, we have the
following obvious equivalences:
$$
\hbox{e is simple} \Leftrightarrow \hbox{ Ker }L_e =\{0\}
\Leftrightarrow  L_e\hbox{ is
surjective }
$$
2) $U(t,0): H^s(S^1) \longrightarrow H^s(S^1)$ is compact and so
$$ e \hbox{ is hyperbolic } \Leftrightarrow \hbox{ no eigenvalue of
}L_e
\hbox{ belongs to the
imaginary axis}.
$$
3) The operator $L_e$ is a sectorial operator with compact resolvent.
Moreover, the sector is locally uniform in the sense that, for any
$f\in\Gm$ and any equilibrium point $e$ of \eqref{eq}, there exist
positive constants $C$ and $C'$ uniform in a neighborhood of $(e,f)$
in $H^s(S^1)\times \Gm$ such that, for any
eigenvalue $\lambda$ of $L_e$,
\begin{equation}\label{eq-lemme-Le-spectre}
\left|\text{Im}(\lambda)\right|\leq C-C'\text{Re}(\lambda)~.
\end{equation}
\end{lemma}

\begin{proof}
1) Since $e\in H^{2}(S^1)$, $\varphi \longmapsto
\varphi+f'_u(x,e,e_x)\varphi +f'_{u_x}(x,e,e_x) \varphi_x$ is bounded
from $H^{2}(S^1)$ into $H^1(S^1)$.  Thus, $L_e$ is a compact
perturbation of $\partial^2_{xx}-Id:H^{2}(S^1)\rightarrow
L^2(S^1)$.  Since $\partial^2_{xx}-Id$ is bijective, it is a
Fredholm operator of index $0$, which implies that $L_e$ is also a
Fredholm operator of index $0$.\\
2) Since $L_e: H^2(S^1) \rightarrow L^2(S^1)$ is a real
operator with compact resolvent, its spectrum consists in a sequence
of eigenvalues of finite multiplicity and its spectrum is symmetric
with respect to the real axis. Due to the smoothing properties of the
parabolic flow, $e^{L_e}$ is
a compact operator, thus its spectrum consists in $\{0\}$ and a
sequence of
eigenvalues with finite multiplicity converging to $0$. Hence,
      due to the spectral mapping theorem (see theorems 2.2.3 and
2.2.4 of \cite{Pazy} for example),
\begin{equation}
\label{spectralmap}
        \exp(t \sigma (L_e)) \subset \sigma (e^{L_e t}) \subset \exp(
t\sigma
        (L_e)) \cup \{0\}~, \quad \forall \,t \ge 0~,
\end{equation}
which clearly implies the equivalence.\\
3) Moreover, $L_e$ is the perturbation of the self-adjoint operator
$\partial_{xx}^2+f'_u(x,e,e_x)$ by the term
$f_{u_x}(x,e,e_x)\partial_x$. Therefore, it is a sectorial operator
by Theorem 3.2.1 of \cite{Pazy} or by Theorem 1.3.2 of \cite{Henry81}.
The fact that the sector is locally uniform in $(e,f)$ and the
estimate \eqref{eq-lemme-Le-spectre} are straightforward consequences
of the proof of Theorem 1.3.2 of \cite{Henry81}. \end{proof}

In Section \ref{subsec-spectre} below, we will prove special spectral
properties of general scalar one-dimensional  linear parabolic
equations. As a direct consequence of
the spectral theorem and of Theorem \ref{A-th-spectre} and
Proposition \ref{A-prop-spectre} below, we obtain the following
interesting spectral property, that we will apply in the proof of
the genericity theorems.

\begin{prop} \label{spectreLe}
Let $e$ be any equilibrium point of \eqref{eq}.
Then, the dimension of the generalized eigenspace
corresponding to the eigenvalues of $L_e$ with same real part is at
most
two.  Moreover, if $\lambda$ is a real eigenvalue of $L_e$, then there
exists no other eigenvalue of $L_e$ with real part equal to $\lambda$.
\end{prop}


\subsection{The case of periodic orbits}

In this section, we consider the non-trivial time-periodic solutions
of Equation \eqref{eq}.  Let $p=p_0(x,t)$ be a periodic solution of
\eqref{eq} of period $T_0 >0$.  We recall that $U(t,0) \equiv
U_f(t,0)$ denotes the linear operator associated with
\eqref{linearized} and that, since $p_0(t)$ is of period $T_0$,
$U(T_0,0)$ is called the {\sl period map}.  \HB
Notice that, since $p_0(t)$ is a periodic
solution with period $T_0$ of the autonomous equation \eqref{eq},
$\mu =1$ is automatically an eigenvalue of $U(T_0,0)$ with
eigenfunction $\partial_t p_0(0)$.

\begin{definition}
Let $p_0(t)$ be a periodic solution of  \eqref{eq}, with period $T_0
>0$. \HB
1) We say that $p_0(t)$ is a simple or non-degenerate periodic orbit
if
the eigenvalue $\mu =1$ is simple and isolated from the rest of the
spectrum. Then, $T_0$ is called a simple period. \HB
2) We say that $p_0(t)$ is hyperbolic if the eigenvalue $\mu=1$ is
simple and if the intersection of the spectrum of $U(T_0,0)$ with the
unit circle reduces to the eigenvalue $1$.
\end{definition}

Remark that, since the equation
\eqref{linearized} is smoothing in finite time, $U(T_0,0)$ is a
compact map from $H^s(S^1)$ into itself.  Therefore,
the spectrum of $U(T_0,0)$ consists of $0$ and a sequence of
eigenvalues $\mu_k$ converging to $0$. Thus, $p_0(t)$ is a hyperbolic
periodic orbit if  $\mu=1$ is the only eigenvalue of $U(T_0,0)$ on
the unit circle and is simple. \HB
Several properties of the spectrum of $U(T_0,0)$
are described in Section 2.3 below. In particular, as a direct
consequence of Theorem \ref{A-th-spectre}, Proposition
\ref{A-prop-spectre} and the fact that $\mu =1$ is an eigenvalue of
$U(T_0,0)$, we obtain the following result.
\begin{prop} \label{simplehyper}
Let $p_0(t)$ be a periodic orbit of \eqref{eq} of period
$T_0$. Then, $\mu =1$ is an eigenvalue of the period map $U(T_0,0)$
of multiplicity at most $2$ and there is no other eigenvalue on the
unit circle. As a consequence, we have the equivalence:
$$
p_0(t) \hbox{ is hyperbolic } \Leftrightarrow  p_0(t) \hbox{ is
simple}.
$$
\end{prop}

\noindent {\bf Nota Bene: } For more general
equations than \eqref{eq}, if for example, $\mu = \exp 2i \pi/n$ is
an eigenvalue of $U(T_0,0)$, the periodic orbit $p_0(t)$, considered
as periodic solution of period $nT_0$, will not be simple. So, in
general, the definition of
simplicity can depend on the chosen period, whereas the definition of
the hyperbolicity does not depend on the chosen period. This will
not be the case here. Indeed, let $p_0(t)$ be a periodic solution of
\eqref{eq} with minimal period
$T_0$. Proposition \ref{simplehyper} implies that the generalized
eigenfunctions associated to the eigenvalue $\mu =1$ of the map
$U(nT_0,0)$
are exactly the eigenfunctions associated to the eigenvalue $\mu =1$
of the
period map $U(T_0,0)$.\\

\subsection{Spectrum of a general linear
operator}\label{subsec-spectre}

In this section, we consider a more general linear equation.
Let $T>0$ be any positive time.  Let $a$ and $b$ be two functions in
$\Cm^1(S^1\times[0,T],\Rm)$.  We consider the operator
$U(T,0): H^s(S^1)\longrightarrow H^s(S^1)$ defined by
$U(T,0)w_0=w(T)$ where $w(t)$ is the solution of
\begin{equation}
\label{A-eq-lin}
\begin{split}
\partial_t w(x,t) &=\partial^2_{xx}w(x,t)
+a(x,t)w(x,t)+b(x,t)\partial_x w(x,t)~, \quad \forall (x,t)\in
S^1\times(0,T] \cr
w(x,0) &=w_0(x)~.
\end{split}
\end{equation}
Due to the classical smoothing properties of parabolic equations,
$U(T,0)$ is a compact operator from $H^s(S^1)$ into itself.  Thus,
the spectrum
consists in $\{0\}$ and a sequence of nonzero eigenvalues of finite
multiplicity converging to $0$.  Notice that $0$ is not an eigenvalue
due to the backward uniqueness property of the parabolic equation.  We
denote by $(\lambda_k)_{k\in\Nm}$ the eigenvalues of $U(T,0)$ with
the convention
that they are repeated according to their multiplicity and ordered by
$|\lambda_{k+1}|\leq|\lambda_k|$. \HB
Using the properties of the zero number, Angenent (\cite{Angenent})
has shown the following result (a first statement
with $a$ and $b$ analytic was proved in \cite{Angenent-Fiedler}).

\begin{theorem}\label{A-th-spectre}
Let $(\lambda_k)_{k\in\Nm}$ be the spectrum of $U(T,0)$ as introduced
above. \HB
Then, for all $j \geq 0$,  $|\lambda_{2j}|>|\lambda_{2j+1}|$.  In
particular, $\lambda_0$ is a simple real eigenvalue.\\
Moreover, let $E_0$ denote the one-dimensional eigenspace
corresponding to $\lambda_0$ and for $j\geq 1$, let $E_{2j}$ denote
the two-dimensional real generalized eigenspace corresponding to
$\{\lambda_{2j-1},\lambda_{2j}\}$.  Then, for all $j\geq 0$, any
nonzero real function $v\in E_{2j}$ has exactly $2j$ zeros and all
these zeros are simple.
\end{theorem}

Assume that $\lambda$ is a real eigenvalue and that $\mu$ is another
eigenvalue such that $|\lambda|=|\mu|$.  Since $a$ and $b$ are real
functions, $\overline \mu$ is also an eigenvalue.  Theorem
\ref{A-th-spectre} shows that there are at most two eigenvalues with
same modulus, thus $\mu$ is also real.  It could be possible that
$\mu=-\lambda$, but the following result prevents this case.

\begin{prop}\label{A-prop-spectre}
Let $j\geq 1$.  If $\lambda_{2j-1}$ and $\lambda_{2j}$ are two
consecutive eigenvalues of $U(T,0)$, then
$\lambda_{2j-1}\lambda_{2j}>0$.  In particular, if $\lambda$ is a real
     eigenvalue then $-\lambda$ is not an eigenvalue
and there is no other distinct eigenvalue on the circle
$\{z\in\Com,~|z|=|\lambda| \}$.
\end{prop}
\begin{proof}
Due to Theorem \ref{A-th-spectre}, there are at most two eigenvalues
with same modulus (counting multiplicity).  Since the spectrum is
symmetric with respect to the imaginary axis, if $\lambda_{2j-1}$ or
$\lambda_{2j}$ is not real, then both must be conjugated and
$\lambda_{2j-1}\lambda_{2j}=|\lambda_{2j}|^2>0$.  For the same reason,
it is also clear that the last assertion of Proposition
\ref{A-prop-spectre} is a direct consequence of the first one.\\
The only case, which has to be studied, is the case when
$\lambda_{2j-1}$ and $\lambda_{2j}$ are
both real.  We argue by contradiction: assume that $\lambda_{2j-1}$
and $\lambda_{2j}$ are real eigenvalues of $U(T,0)$ of opposite signs.
We denote by $U_\varepsilon(T,0)$ the evolution operator corresponding
to
\begin{equation}\label{A-eq-lin-eps}
\partial_t w(x,t)=\partial^2_{xx}w(x,t)+\varepsilon a(x,t)w(x,t)
+\varepsilon b(x,t)
\partial_x w(x,t)~, \quad \forall (x,t)\in S^1\times(0,T]~.
\end{equation}
We denote by $(\lambda^\varepsilon_k)_{k\in\Nm}$ the set of
eigenvalues of $U_\varepsilon(T,0)$.  We set
$$\Em=\{\varepsilon \in
[0,1] \, | \, \lambda^\varepsilon_{2j-1}\text{ and
}\lambda^\varepsilon_{2j}\text{ are real eigenvalues of
}U_\varepsilon(T,0)\text{ of opposite signs}\}~.
$$
By assumption, $1\in\Em$ and trivially $0\not\in\Em$.  We will obtain
a contradiction by showing that $\Em$ is open and closed.\\
\underline {O}p\underline{enness:} let $\varepsilon_0\in\Em$.  Since
$\lambda_{2j}^{\varepsilon_0}$ and $\lambda_{2j-1}^{\varepsilon_0}$
are different,
they must be simple
due to Theorem \ref{A-th-spectre}.  Let $\varphi_{2j}^{\varepsilon_0}$ and
$\varphi_{2j-1}^{\varepsilon_0}$ be two corresponding normalized real
eigenfunctions.
They have $2j$ simple zeros.  The simplicity of the eigenvalues and
the implicit function theorem imply that there exists a neighborhood
$\Vm$ of $\varepsilon_0$ and functions $\mu_1(\varepsilon)$,
$\mu_2(\varepsilon)$, $\psi_1(\varepsilon)$ and
$\psi_2(\varepsilon)$ defined on $\Vm$ such that the following
properties hold.  The functions $\mu_i$ are of
class $\Cm^0(\Vm,\Rm)$, $\mu_1(\varepsilon_0)=\lambda_{2j-1}^{\varepsilon_0}$,
$\mu_2(\varepsilon_0)=\lambda_{2j}^{\varepsilon_0}$ and
$\mu_i(\varepsilon)$ is a
simple real eigenvalue of $U_\varepsilon(T,0)$.  The functions
$\psi_i$ are of class $\Cm^0(\Vm,H^s(S^1))$,
$\psi_1(\varepsilon_0)=\varphi_{2j-1}^{\varepsilon_0}$,
$\psi_2(\varepsilon_0)=\varphi_{2j}^{\varepsilon_0}$ and
$\psi_i(\varepsilon)$ is the
normalized real eigenfunction corresponding to $\mu_i(\varepsilon)$.
Since the zeros of all eigenfunctions of $U_\varepsilon(T,0)$ are simple and
$H^s(S^1)\subset\Cm^1(S^1)$, restricting $\Vm$ if necessary, we can
assume that $\psi_i(\varepsilon)$ has exactly $2j$ zeros which are all
simple.  Due to Theorem \ref{A-th-spectre}, this shows that for all
$\varepsilon\in \Vm$, the eigenvalues $\mu_i(\varepsilon)$ are the
eigenvalues $\lambda^\varepsilon_{2j-1}$ and
$\lambda^\varepsilon_{2j}$ of $U_\varepsilon(T,0)$.  Moreover,
the eigenvalues of $U_\varepsilon(T,0)$
     cannot be zero due to the backward uniqueness property of
\eqref{A-eq-lin-eps}.  Thus, up to a restriction of $\Vm$, the
eigenvalues $\mu_i(\varepsilon)$ are of opposite signs and
$\Vm\subset\Em$.

\noindent \underline{Closeness:} assume that
$(\varepsilon_n)\subset\Em$ is a
given sequence such that
$\varepsilon_n\longrightarrow \varepsilon$.  We first consider the
sequence of eigenvalues $\lambda^{\varepsilon_n}_{2j}$ and the
sequence of corresponding real eigenfunctions
$\varphi^{\varepsilon_n}_{2j}$ with
$\|\varphi^{\varepsilon_n}_{2j}\|_{L^2(S^1)}=1$.  Taking the inner
product in $L^2$ of \eqref{A-eq-lin-eps} with $w$, integrating in time
and applying Gronwall Lemma, we show that
there exists a positive constant $C_0$ independent of $\varepsilon$ such that
every solution $w$ of \eqref{A-eq-lin-eps} satisfies
$$
\|w(T)\|_{L^2(S^1)}\leq C_0\|w(0)\|_{L^2(S^1)}~.
$$
Thus, $|\lambda^{\varepsilon_n}_{2j}|$ is bounded by $C_0$ and, up to
the
extraction of a subsequence, we can assume that
$\lambda^{\varepsilon_n}_{2j}$ converges to $\mu$.  First, $\mu$
cannot be zero.  Indeed, there exists a closed curve
$\Gamma\subset\Com \setminus\{0\}$ surrounding the $2j+2$ first
eigenvalues of $U_\varepsilon(T,0)$ but not zero.  The projector
$$
P_{2j+2}= -\frac{1}{2i\pi}\int_\Gamma (U_\varepsilon(T,0)-z)^{-1}dz
$$
is the spectral projector of $U_\varepsilon(T,0)$ onto the space
generated by the first $2j+2$ eigenvalues.
Thus, for $n$ large enough
$$
P^n_{2j+2}= -\frac{1}{2i\pi}\int_\Gamma
(U_{\varepsilon_n}(T,0)-z)^{-1}dz
$$
is a spectral projector of rank
$2j+2$ and thus $\Gamma$ still surrounds $2j+2$ eigenvalues of
$U_{\varepsilon_n}(T,0)$.  Due to the ordering $|\lambda_0|\geq
|\lambda_1|\geq ...$, this shows that $\lambda^{\varepsilon_n}_{2j}$
must be away from zero for large $n$.  Next, using the
smoothing property of parabolic equations shows that there exists
a positive constant constant $C_1$ independent of $\varepsilon$ such
that every solution $w$
of \eqref{A-eq-lin-eps} satisfies
$$
\|w(T)\|_{H^2(S^1)}\leq
C_1\|w(0)\|_{L^2(S^1)}~.
$$
Therefore,
$$
|\lambda^{\varepsilon_n}_{2j}|\|\varphi^{\varepsilon_n}_{2j}\|_{H^2(S^1)}\leq
C_1\|\varphi^{\varepsilon_n}_{2j}\|_{L^2(S^1)}=C_1~.
$$
Since
$\lambda^{\varepsilon_n}_{2j}$ is bounded away from zero,
$(\varphi^{\varepsilon_n}_{2j})$ is bounded in $H^2(S^1)$ and we can
assume that $(\varphi^{\varepsilon_n}_{2j})$ converges to a function
$\psi$ in $\Cm^1(S^1)$.  Passing to the limit, we see that $\mu$ is a
real eigenvalue of $U_\varepsilon(T,0)$ with normalized eigenfunction
$\psi$.  Due to Theorem \ref{A-th-spectre}, $\psi$ has simple zeros
only and since $(\varphi^{\varepsilon_n}_{2j})$ converges to $\psi$ in
$\Cm^1(S^1)$, the number of zeros of $\psi$ must be $2j$.  Hence,
$\mu$
must be equal either to $\lambda^\varepsilon_{2j}$ or to
$\lambda^\varepsilon_{2j-1}$.
Arguing in the same way for $\lambda^{\varepsilon_n}_{2j-1}$ we show
that, up to the extraction of a subsequence,
$\lambda^{\varepsilon_n}_{2j}$ and $\lambda^{\varepsilon_n}_{2j-1}$
converge to $\lambda^{\varepsilon}_{2j}$ or
$\lambda^{\varepsilon}_{2j-1}$.  Since neither
$\lambda^{\varepsilon}_{2j}$ nor $\lambda^{\varepsilon}_{2j-1}$ can be
zero and since
$\lambda^{\varepsilon_n}_{2j}\lambda^{\varepsilon_n}_{2j-1}<0$, both
$\lambda^{\varepsilon}_{2j}$ and $\lambda^{\varepsilon}_{2j-1}$ are
limits of one of the sequences $(\lambda^{\varepsilon_n}_{2j})$ or
$(\lambda^{\varepsilon_n}_{2j-1})$ and thus are real and of opposite
signs.  Therefore, $\varepsilon$ belongs to $\Em$.\\
All these properties yield a contradiction since $[0,1]$ is
connected. The proposition is proved. \end{proof}


\section{Generic hyperbolicity of the equilibrium
points}\label{Sec-eq}

\subsection{Generic simplicity}

Generic simplicity of equilibria has been proved for various
parabolic systems, including the Navier-Stokes equations,
in different settings (see \cite{Quinn}, \cite{sauttemam}, \cite{SW},
\cite{Babin-V-livre}, \cite{Bruno-Chow}, \cite{Bruno-Pola} for
example). Since the proof is rather simple, we include it here for
sake of completeneness. We follow the lines of
the proof given by  \cite{Bruno-Pola} in the case where $f$ does not
depend on the derivative $u_x$.

\begin{prop}\label{prop-gen-simple}
For any $n\geq 1$, the set
$$
\Om^{\text{simp}}_n=\{f\in\Gm \, | \, \text{ any
equilibrium point }e\text{ of \eqref{eq} with }\|e\|_{\Cm^1(S^1)}\leq
n\text{ is simple}\}
$$
is a dense open subset of $\Gm$.  As a consequence, the set
$$
\Om^{\text{simp}}=\{f\in\Gm \, | \, \text{ any equilibrium point }e\text{ of
\eqref{eq} is simple}\}
$$
is a generic subset of $\Gm$.
\end{prop}

\begin{proof}
If $\Om^{\text{simp}}_n$ is a dense open set, then, since any equilibrium
belongs to $\Cm^1(S^1)$, $\Om^{\text{simp}}=\cap_n \Om^{\text{simp}}_n$ and
thus $\Om^{\text{simp}}$ is a generic subset of $\Gm$.\\
Let $n$ be any positive integer.\\
$\underline\Om^{\text{simp}}_n$\underline{ is o}p\underline{en:}
  We denote by $R: g \in \Gm \mapsto Rg \in C^2(S^1 \times [-(n+2), n
  +2] \times [-(n +2),
n+2],\Rm)$ the restriction operator defined by
$$
Rg =  g_{| S^1 \times [-(n+2), n+2] \times [-(n +2), n+2]}~.
$$
We notice that $R$ is a continuous, surjective map and that, on
$R\Gm=C^2(S^1 \times [-(n +2), n +2] \times [-(n +2), n +2],\Rm)$,
the Whitney topology and the
classical $C^2$-topology coincide. Since $\Om^{\text{simp}}_n$ only
depends on the values of $f$ in $S^1 \times [-n,n] \times [-n,n]$, it
suffices to show that $R\Om^{\text{simp}}_n$ is open in $R\Gm$ for
the classical $C^2$-topology. It thus suffices to prove, that if
$(f_k)$ is a sequence of functions in
$\Gm\setminus\Om^{\text{simp}}_n$ converging to $f$ in $C^2(S^1
\times [-(n +2), n+2] \times [-(n +2), n +2],\Rm)$,
then $f$ belongs
to $\Gm\setminus\Om^{\text{simp}}_n$.
Assume that $(f_k)$ is such a
sequence of functions in
$\Gm\setminus\Om^{\text{simp}}_n$ converging to $f$.
Due to Lemma
\ref{lemme-Le-Fredholm}, there exist a sequence of equilibria $(e_k)$
with $\|e_k\|_{\Cm^1}\leq n$ and a sequence of functions
$(\varphi_k)\subset H^2(S^1)$ with $\|\varphi_k\|_{L^2}=1$ and
$L_{e_k}\varphi_k=0$.  Since
$$
\int |\partial_x\varphi_k|^2=\int
D_uf_k(x,e_k,\partial_x e_k)|\varphi_k|^2+ \int
D_{u_x}f_k(x,e_k,\partial_x
e_k)\partial_x\varphi_k\overline\varphi_k~,
$$
the sequence $(\varphi_k)$ is bounded in $H^1(S^1)$.  Moreover,
$\partial^2_{xx}e_k=f_k(x,e_k,\partial_x e_k)$ and thus the sequence
$(e_k)$ is bounded in $H^{s+2}(S^1)$.  Up to the extraction of a
subsequence, we can assume that $(e_k)$ converges in $H^2(S^1)$ to
$e$ and $(\varphi_k)$ converges in $\Lm^2(S^1)$ to $\varphi$.  Passing
to the limit, we obtain $L_e\varphi=0$ and $\|e\|_{\Cm^1}\leq n$.
By elliptic estimates, $\varphi\in H^2(S^1)$, thus
$f\not\in\Om^{\text{simp}}_n$, and so $\Om^{\text{simp}}_n$ is open.\\
$\underline\Om^{\text{simp}}_n$\underline{ is dense:} let $f\in\Gm$ be
given.  We follow the now classical method (see \cite{Bruno-Pola} for
example).
We introduce a function $\chi\in\Cm^2(\Rm^2,\Rm)$ satisfying
$\chi\equiv 1$ in $[-n-1,n+1]^2$ and $\chi\equiv 0$ in
$\Rm\setminus[-n-2,n+2]^2$.  For any open neighborhood $\Vm$ of $f$ in
$\Gm$, there exists an open neighborhood $\Um$ of $0$ in $\Cm^2(S^1)$
such
that, for all $a\in\Um$, the function $(x,u,v)\mapsto
f(x,u,v)+a(x)\chi(u,v)$ belongs to $\Vm$.  So, it is sufficient to
prove that there is a dense set of functions $a\in\Cm^2(S^1,\Rm)$ such
that $f+a\chi\in \Om^{\text{simp}}_n$.  To obtain this density, we apply
Sard-Smale Theorem (Theorem \ref{A-th-Sard-Smale} in the appendix),
to the functional $\Phi : (e,a) \in \{e\in
H^{2}(S^1),\|e\|_{\Cm^1}<n+1\}\times\Cm^2(S^1,\Rm) \mapsto \Phi(e,a)
\in L^2(S^1)$ defined by
$$
\Phi(e,a) = e_{xx}+f(x,e,e_x)+a(x)\chi(e,e_x)
$$
and to the point $z=0$. \HB
Assumption iii) of Theorem \ref{A-th-Sard-Smale}
is satisfied.  To simplify the notation, let us prove i) and ii) with
$a=0$, which does not lead to any loss of generality.  First notice
that $(e,0)\in\Phi^{-1}(0)$ means that $e$ is an equilibrium point of
\eqref{eq}.  Assumption i) of Theorem \ref{A-th-Sard-Smale} is
satisfied since $D_e\Phi(e,0)=L_e$ is a Fredholm operator of
index $0$ due to Lemma \ref{lemme-Le-Fredholm}.  To prove ii), we must
find for each $h\in L^2(S^1)$ a pair $(\varphi,b)$ such that
$D_{e,a}\Phi(e,0).(\varphi,b)=L_e\varphi+b(x)\chi(e,e_x)=L_e\varphi+b(x)=h$.
Since $L_e$ is Fredholm, $h-b\in Im(L_e)$ if and only if, for all
$\psi\in Ker(L_e^*)$, $\int(h-b)\psi=0$.  As $Ker(L_e^*)$ is
finite-dimensional, we can introduce a finite orthonormal basis
$(\psi_i)_{i=1\ldots p}$ of $Ker(L_e^*)$.  We are reduced to find $b$
such that $\int(h-b)\psi_i=0$ for all $i=1,\ldots,p$ that is, so that $0$
belongs to the image of the map: $b\in\Cm^2(S^1)\longmapsto
(\int(h-b)\psi_i)_{i=1\ldots p}$.  This image is closed since it is an
affine subspace of $\Rm^p$ and by density of $\Cm^2$ in $L^2$, we
can find $b$ as close to $h$ as wanted.  Therefore, we can find $b$
such that $\int(h-b)\psi_i=0$ for all $i=1,\ldots,p$ and thus
$D_{e,a}\Phi$ is surjective. Actually, remarking that
the vectors $\psi_i$ belong to $H^3(S^1)$, we can simply choose $b=
\sum_{0}^{p}(h,\psi_i) \psi_i$. \\
The conclusion of Theorem \ref{A-th-Sard-Smale} shows that for a
generic $a\in\Cm^2(S^1)$, any $e$ satisfying $\|e\|_{\Cm^1}\leq n$ and
$e_{xx}+f(x,e,e_x)+a(x)\chi(e,e_x)=0$ is such that $D_e\Phi=L_e$ is
surjective which implies by Lemma \ref{lemme-Le-Fredholm} that $e$ is
a simple equilibrium point.  Thus, we can choose a function $a$ as
small as wanted such that $f+a\chi\in \Om^{\text{simp}}_n$, which shows
that $\Om^{\text{simp}}_n$ is dense.
\end{proof}

\subsection{Generic hyperbolicity}

In Proposition \ref{prop-gen-simple}, we have proved the generic
simplicity
of the equilibrium points of \eqref{eq}.  Using this result, it is not
difficult to show the generic hyperbolicity.

\begin{prop}\label{prop-gen-hyper}
For any $n\geq 1$, the set
$$
\Om^{\text{h}}_n=\{f\in\Gm\, | \,\text{ any
equilibrium point }e\text{ of \eqref{eq} with }\|e\|_{\Cm^1(S^1)}\leq
n\text{ is hyperbolic}\}
$$
is a dense open subset of $\Gm$.  As a consequence, the set
$$
\Om^{\text{h}}=\{f\in\Gm \, | \,\text{ any equilibrium point }e\text{ of
\eqref{eq} is hyperbolic}\}
$$
is a generic subset of $\Gm$.
\end{prop}

\begin{proof}
As in the proof of Proposition \ref{prop-gen-simple}, it is sufficient
to show that $\Om^{\text{h}}_n$ is a dense open set.\\

$\underline\Om^{\text{h}}_n$\underline{ is o}p\underline{en:}
As in the proof of Proposition \ref{prop-gen-simple}, it is
sufficient to prove that if
$(f_k)$ is a sequence of functions in
$\Gm\setminus\Om^{\text{h}}_n$ converging to $f$ in $C^2(S^1 \times [-(n+2),
n+2] \times [-(n+2), n+2],\Rm)$,
then $f$ belongs
to $\Gm\setminus\Om^{\text{h}}_n$.
So let $(f_k)$ be such a
sequence of functions in
$\Gm\setminus\Om^{\text{h}}_n$ converging to $f$.
Then there
exist sequences
$(e_k)\subset H^{s+2}(S^1)$, $(\lambda_k)\subset i\Rm$ and
$(\varphi_k)\subset H^{2}(S^1)$ such that
$\partial^2_{xx}e_k+f_k(x,e_k,\partial_x e_k)=0$, $\|e_k\|_{\Cm^1}\leq
n$, $\|\varphi_k\|_{L^2}=1$ and
$L_{e_k}\varphi_k=\lambda_k\varphi_k$.  As in the proof of Proposition
\ref{prop-gen-simple}, we can assume that $(e_k)$ converges in
$H^2(S^1)$ to a function $e$.  Then, the proof of openness will be
exactly the same as the corresponding one in Proposition
\ref{prop-gen-simple} as
soon as we show that there exists a subsequence
     $(\lambda_{k_n})\subset i\Rm$, which  is convergent.  This is
a consequence of Estimate \eqref{eq-lemme-Le-spectre}.\\
$\underline\Om^{\text{h}}_n$\underline{ is dense:} let $f_0\in\Gm$.
We are going to show that we can construct successive perturbations
of $f_0$, as
small as needed, to obtain a function $f\in\Om^{\text{h}}_n$ as
closed to $f_0$ as is wanted.  First, due to Proposition
\ref{prop-gen-simple}, we can find $f_1$ close to $f_0$ such that all
equilibrium points of $f_1$ are simple.  The set
$\{e \, | \, \partial^2_{xx}e+f_1(x,e,e_x)=0,~ \|e\|_{\Cm^1}\leq
n+1\}$ is
bounded in $H^{s+2}(S^1)$ and hence compact in $H^2(S^1)$.  Since all
equilibrium points of $f_1$ are simple and thus isolated, there is a
finite number of equilibria $e_1,\ldots,e_p$ of $f_1$ which satisfy
$\|e_j\|_{\Cm^1}\leq n+1$.\\
We next explain how each equilibrium $e_i$ can be made hyperbolic
by successive perturbations of $f_1$. Let $e$ be an equilibrium
point of \eqref{eq} with $\|e\|_{\Cm^1}\leq n$ and let us denote
$(\lambda_k)_{k\in\Nm}$ the sequence of eigenvalues of the
corresponding linearized operator $L_e$. Let
$\chi\in\Cm^2(\Rm,\Rm)$ be a smooth cut-off function such that, for
example, $\chi(y)=y$ for $|y| \leq n+1$ and $\chi(y)=0$ for $|y| \geq
n+2$. Then, $\chi(e(x))=e(x)$ for any $x \in S^1$.
If we perturb $f_1$ by setting for small
$\alpha\in\Rm$, $f_\alpha(x,v,w)=f_1(x,v,w)+\alpha(\chi(v)-e(x))$, then
$e$ is still an equilibrium point of \eqref{eq} with $f$ replaced by
$f_\alpha$ and the spectrum of $L_e$ becomes
$(\lambda_k+\alpha)_{k\in\Nm}$.  Since
the eigenvalues of $L_e$ are isolated, this means that one can
perturb $f_1$ such that
$L_e$ has no longer eigenvalues on the imaginary axis, that is, such
that $e$ becomes hyperbolic. On the other hand, by the implicit functions
theorem, if $e$ is a simple (resp.  a hyperbolic) equilibrium of $f$,
there exist neighborhoods $\Vm\subset\Cm^1(S^1)$ of $e$ and
$\Um\subset\Gm$ of $f$ such that for all $g\in\Um$, there exists a
unique equilibrium point $e(g)\in\Vm$ and this equilibrium is simple
(resp. hyperbolic). Thus, we can make successively each equilibrium
hyperbolic without changing the status of the other equilibria
$e_1,\ldots,e_p$.
\end{proof}


\section{Generic hyperbolicity of the periodic orbits}\label{Sec-per}

The aim of this section is the proof of the genericity theorem
\ref{Peixoto}, that is, we want to show that there exists a subset
$\Om$
of $\Gm$ such that, for all $f \in \Om$, all the equilibrium points
and all the periodic solutions of \eqref{eq} are hyperbolic.

To show Theorem \ref{Peixoto}, we use the induction argument of
Peixoto \cite{Pei} in the case of vector fields on compact manifolds
or of Mallet-Paret \cite{Mallet} in the case of functional
differential
equations. Thus, for any $n \geq 1$ and any $A>0$, we introduce the
set
\begin{equation*}
\begin{split}
\Om(A,n) = &\{f \in \Om^{h}_{n} \, | \, \hbox{ all nonconstant
periodic
orbits } p(t) \hbox{ of \eqref{eq} with period } \cr
& T \in (0,A]
\hbox{ such that sup}_{t\in\Rm} \|p(t)\|_{C^1(S^1)} \leq n \hbox{ are
non-degenerate} \}~.
\end{split}
\end{equation*}
Due to Proposition \ref{simplehyper}, we have the following equality:
\begin{equation*}
\begin{split}
\Om(A,n) =&\{f \in \Om^{h}_{n} \, | \, \hbox{ all nonconstant periodic
orbits } p(t) \hbox{ of \eqref{eq} with period } \cr
& T \in (0,A]
\hbox{ such that sup}_{t\in\Rm} \|p(t)\|_{C^1(S^1)} \leq n \hbox{ are
hyperbolic} \}.
\end{split}
\end{equation*}
As in \cite{Pei} and in \cite{Mallet}, we
show that $\Om(n,n)$ is open and dense in $\Gm$. Since $\Gm$ is a
Baire space,
it follows that $\Om = \cap_{n=1}^{\infty} \Om(n,n)$ is a
generic subset and hence a dense subset of $\Gm$, which proves
Theorem \ref{Peixoto}.
\vskip 2mm

Let $n\in\Nm$ be fixed for the remaining part of the section.  We
first begin with auxiliary results, which will be widely used in this
section.  In particular, we prove
that, for any $A>0$, $\Om(A,n)$ is open. To simplify our statements
below, we denote by $S_{f}(t)$ the local nonlinear semigroup defined
by the equation \eqref{eq} with nonlinearity $f$.

\begin{prop} \label{period-mini}
The following properties hold. \HB
(a) Let $f\in\ \Om^{\text{simp}}_n$ (resp.
$f\in\Om^h_n$). There exists $\delta>0$ such that
$f\in\Om^{\text{simp}}_{n+\delta}$ (resp. $f\in\Om^h_{n+\delta}$).\HB
(b) For any $\mu>0$, the set $\Om(A,\mu)$ is open. \HB
(c) Let $f \in \Om^{\text{simp}}_n$ and let $\delta$ be
as in (a). There exist $\varepsilon>0$ and a neighborhood $\Nc_1 \subset
\Om^{\text{simp}}_{n+\delta}$ of $f$ such that, for any $g \in \Nc_1$, any
nonconstant periodic solution
$p(t)$ of $S_{g}(t)$ with $\sup_{t \in \Rm} \|p(t)\|_{C^1(S^1)} \leq
n +\delta$ has a
smallest period strictly larger than $\varepsilon$. \HB
(d)  Let $f \in \Om_{n}^{h}$ and let $\delta$ be
as in (a). For any $A>0$, there exist a positive constant $r$
and a neighborhood  $\Nc_2 \subset
\Om^{h}_{n+\delta}$  of $f$ such that the following property holds.
Let $\Em_{f,n+\delta}$ be the set of the equilibrium points $e_f$ of $S_f(t)$
satisfying $\|e_f\|_{C^1(S^1)}\leq n+\delta$ and
let $\Bc\Em(f,n+\delta,r)=\cup_{e_f\in\Em_{f,n+\delta}} B_{H^s(S^1)}(e_f,r)$.
For any $g \in \Nc_2$, $\Em_{g,n+\delta}\subset\Bc\Em(f,n+\delta,r)$ 
and the set
of all nonconstant periodic orbits $p(t)$ of $S_g(t)$ of period less 
than $A$ and
satisfying $\sup_{t \in\Rm} \|p(t)\|_{C^1(S^1)} \leq n+\delta$ does not
intersect $\Bc\Em(f,n+\delta,r)$.
\end{prop}

\begin{proof}
The proof of property (a) is similar to the proof of the openness of
$\Om^{\text{simp}}_n$ in Proposition \ref{prop-gen-simple}.

To prove property (b), we follow the lines of the proof of Theorem 2.1 of
\cite{Mallet}.
By the remarks made at the beginning of the proof of Proposition 
\ref{prop-gen-simple}, it is
sufficient to prove that, if $f_m$ is a sequence of functions in
  $\Gm \setminus \Om(A,\mu)$ converging to $f$ in $C^2(S^1 \times 
[-(\mu +2), \mu
  +2] \times [-(\mu +2), \mu +2],\Rm)$,
then $f$ belongs to $\Gm \setminus \Om(A,\mu)$. So assume that
$f_m$ is such a sequence of functions in $\Gm
\setminus \Om(A,\mu)$ converging to $f$ in $C^2(S^1 \times[-(\mu +2), \mu +2]
\times [-(\mu +2), \mu +2],\Rm)$.
Then, for any $m$,
there exists a non-simple periodic solution $p_m(t)$ of
\begin{equation}
\label{eqm}
u_t(x,t) =u_{xx}(x,t)+f_m(x,u,u_x)~, \quad (x,t)\in S^1\times\Rm_+^*~,
      \end{equation}
with non-simple period $T_m \in (0, A)$. Moreover, $p_m(t)$ satisfies
the following bound
\begin{equation}
\label{bornen}
\sup_{t} \|p_m(t)\|_{C^1(S^1)} \leq \mu~.
\end{equation}
Due to the smoothness hypotheses made on $f$ and to the bound
\eqref{bornen}, we can show by a recursion argument that the function $p_m(t)$
belongs to $C^1([0,2A], H^s(S^1)) \cap C^0([0,2A], H^{2}(S^1))$ and
that there exists a positive constant $C_0$ such that , for any $m
\in \Nm$,
\begin{equation}
\label{borneH4s+2}
\|p_m(t)\|_{C^1([0,2A], H^s(S^1)) \cap C^0([0, 2A], H^{2}(S^1))}
\leq
C_0~.
\end{equation}
Thus, the family of mappings $p_m(t) \in C^1([0, 2A], H^s(S^1))$, $m
\in
\Nm$ is equicontinuous from $[0,A]$ into $H^s(S^1)$.   By the
Ascoli theorem, there exists a subsequence $(p_{m_j}(t))$ which
converges to a function $p(t)$ in $C^0([0, 2A], H^s(S^1))$
and $p(t)$ belongs to $C^0([0, 2A], H^s(S^1))$.
Furthermore, $f_{m_j}(x,p_{m_j}(t), \partial_x p_{m_j}(t))$ converges to
$f(x,p(t), \partial_x p(t))$ in $C^0([0, 2A], L^2(S^1))$. Taking now
the
limits in the variation of constants formula
$$
p_{m_j}(t) \, = \, e^{B t}p_{m_j}(0) + \int_{0}^{t}
e^{B(t-s)}(p_{m_j}(s) + f_{m_j}(x, p_{m_j}(s), \partial_x
p_{m_j}(s)) )ds ~,
$$
where $B$ is the selfadjoint operator $B =\partial_{xx} -I$, we
conclude that $p(t)$ is a mild solution
of \eqref{eq} on the time interval $[0,2A]$ that is, of the equation
\begin{equation}
\label{mild}
p(t) \, = \, e^{B t}p(0) + \int_{0}^{t}
e^{B(t-s)}(p(s) + f(x, p(s), \partial_x
p(s)) )ds ~.
\end{equation}
Since $p_m(t)$ is a periodic solution of period $T_m \leq A$, we also
conclude that $p(t)$ is a periodic solution of \eqref{eq} of period
$T^*= \lim_{m \rightarrow +\infty} T_m$ and the integral equality
\eqref{mild} holds on $\Rm^+$. Furthermore $p(t)$ is a classical
solution of \eqref{mild}. \HB
We next consider the linearized equations
\begin{equation}
\label{eqlineairem}
\begin{split}
\varphi_t(x,t) &=\varphi_{xx}(x,t)+
D_{u}f_m(x,p_m,\partial_xp_m)\varphi
+ D_{u_x}f_m(x,p_m,\partial_xp_m)\varphi_x~, \cr
\varphi (x,0) &= \varphi_0(x)~.
\end{split}
\end{equation}
as well as the associated linear operators $U_m(t,0)$ defined by $U_m(t,0)
\varphi_0=
\varphi^m(t)$ where $\varphi^m(t)$ is the solution of
\eqref{eqlineairem}. We recall that $U_m(T_m,0) : H^s(S^1) \mapsto
H^s(S^1)$ is a compact map. By assumption, the element $1$, which
belongs to the spectrum $\sigma(U_m(T_m,0))$ of $U_m(T_m,0)$, is of
multiplicity greater than one. \HB
We also consider the operator $U(t,0)$ associated to the limiting
linearized equation
\begin{equation*}
\begin{split}
\varphi_t(x,t) &=\varphi_{xx}(x,t)+ D_{u}f(x,p,\partial_xp)\varphi
+ D_{u_x}f(x,p,\partial_xp)\varphi_x~, \cr
\varphi (x,0) &= \varphi_0(x)~.
\end{split}
\end{equation*}
Since $p_{m_j}(t)$ converges to a function $p(t)$ in $C^0([0, 2A],
H^s(S^1))$ and that $f_m$ converges to $f$ in
$C^2(S^1\times[-(\mu +2), \mu +2]
\times [-(\mu +2), \mu +2],\Rm)$, we easily prove
that the operator $U_{m_j}(t,0)$ converges to $U(t,0)$ in
$\Lc(H^s(S^1),H^s(S^1))$, uniformly in $t \in [0,2A]$.  Assume that $T^*
\ne 0$.  Then, $U_{m_j}(T_{m_j},0)$ converges to $U(T^*,0)$ in
$\Lc(H^s(S^1),H^s(S^1))$. \HB
Now three cases may arise: \HB
1) the periodic orbit $p(t)$ is nonconstant and $T^* >0$. The element
$1 \in \sigma(U(T_m,0))$ has multiplicity greater than one, so $1 \in
\sigma(U(T^*,0))$ also
has  multiplicity greater than one. Thus $p(t)$ is a nonsimple
periodic solution with the nonsimple period $T^*$.  And $f  \notin
\Om(A,\mu)$. \HB
2) The periodic orbit $p(t)$ is a constant $a$ and $T^* >0$. The
element
$1$ belongs to the spectrum $\sigma(U(T_m,0))$ of $U(T_m,0)$. Thus,
$1$ also belongs to the spectrum $\sigma(U(T^*, 0))$ of $U(T^*, 0)$.
This implies that $a$ is not hyperbolic. Thus, $f \notin \Om^{h}_{\mu}$
and $f \notin \Om(A,\mu)$. \HB
3)  The periodic orbit $p(t)$ is a constant $a$ and $T^*=0$. For any
$T_0 \in (0,2A]$, we can find a
sequence of positive integers $N_m$ such that $N_mT_m$ converges to
$T_0 $.  Then, $U_{m_j}(N_{m_j}T_{m_j},0)$ converges to $U(T_0,0)$ in
$\Lc(H^s(S^1),H^s(S^1))$. Thus, $1$ belongs to the spectrum
$\sigma(U(T_0,0))$ of $U(T_0,0)$. Since, as explained in Section 2.2,
$U(T_0,0) = e^{L_a T_0}$ and that $L_a$ is a sectorial operator, the
inclusions \eqref{spectralmap} hold, with $e$ replaced by $a$ and
with $t$ replaced by any $T_0 \in (0,2A]$. It follows that $0$
belongs to the spectrum of $L_a$, which means that $a$ is not a
simple equilibrium point and that $f$ does not belong to the set
$\Om^{\text{simp}}_{\mu}$ and thus does not belong to $\Om(A,\mu)$.

To prove the assertions (c) and (d), we argue as in the
above cases 2) and 3). Proposition 4.1 is proved.
\end{proof}

Before entering into the proof of the density of the set
$\Om(n,n)$, $A \leq n$, we want to emphasize that the two main
ingredients of this proof are the Sard-Smale theorem and the
properties of the zero number. The properties of the zero number are
recalled in the introduction (for more details, see
\cite{Sturm, Angenent, Matano82} for example). These properties are used here
through Proposition \ref{simplehyper} and the following primordial
lemma.
\begin{lemma} \label{injectif}
Let $p_0(t)$ be a periodic solution of \eqref{eq} of minimal period
$T_0 >0$.  The map
$$
(x,t) \in S^1 \times [0, T_0) \mapsto (x, p_0(x,t), \partial_x
p_0(x,t))
$$
is one to one.
\end{lemma}

\begin{proof} Assume that this map is not injective. Then there exist
$x_0$, $t_0 \in [0, T_0)$ and $t_1 \in [0, T_0)$, $t_0 \ne t_1$ such
that
$$
p_0(x_0, t_0)=p_0(x_0, t_1)~, \quad \partial_x p_0(x_0, t_0)=
\partial_x p_0(x_0, t_1)~.
$$
The function $v(x,t)= p_0(x, t + t_1 -t_0)-p_0(x,t)$ is a solution of
the equation
$$
v_t (x,t)= v_{xx}(x,t) + b(x,t)v(x,t) + d(x,t) v_x(x,t)~,
$$
where $b(x,t)= \int_{0}^{1} D_uf(x, p_0(t) + s(p_0(t +
t_1-t_0)-p_0(t)), \partial_x p_0(t + t_1-t_0)) ds$ and
$d(x,t) =\int_{0}^{1} D_{u_x}f(x, p_0(t), \partial_x (p_0(t) +s(p_0(t
+ t_1-t_0)
-p_0(t)))) ds$. Moreover, the function $v(x,t)$ satisfies
$v(x_0,t_0)=0$ and $\partial_x v(x_0, t_0)=0$ and does not vanish
everywhere since $|t_1-t_0|<T_0$. Thus,
the zero number $z(v(t))$ drops strictly at $t=t_0$. Since $v(t)$ is
a periodic function of period $T_0$, this leads to a contradiction
with the fact that $z(v(t+T_0))=z(v(t))$. The lemma is proved.
\end{proof}

To prove that $\Om(n,n)$ is dense in $\Gm$, we argue as follows. Let
$f_0\in\Om^h_n$. Due to Proposition \ref{period-mini}, there exists a
constant $\delta\in (0,1)$ such that $f_0 \in \Om^{h}_{n+\delta}$. 
Let $\delta\in(0,1)$,
$\varepsilon>0$ and $\Nc_1$ be as in Proposition \ref{period-mini} (a)
and (c). And for $A=n+2$, let $r>0$, $\Bc\Em(f_0,n+\delta,r)$ and 
$\Nc_2$ be as in Proposition
\ref{period-mini} (d). We next choose a neighborhood $\Nc^0= \Nc_1
\cap \Nc_2$ of $f_0$. Since $\Om^{h}_{n+\delta}$ is
open, we can also assume that $\Nc^0\subset \Om^h_{ n+\delta}$. Let
$\eta$ be a small
positive constant such that $\eta \sum_{j=1}^{+\infty} j^{-4} \leq 
\delta/2$. We set
$$
\Pc(k)=\Om\big((3/2)^k \varepsilon~,~ n+\delta -\eta \sum_{j=1}^{k}
j^{-4}\big)~.
$$
Obviously, one can choose $\varepsilon$ as small as needed so that there
exists $k_0 \in \Nm$ such that  $n \leq (3/2)^{k_0} \varepsilon <n+2$.
In Proposition \ref{generic3A/2A} below, we show that for all $k\in
\Nm$, $0 < k \leq k_0$, there exists an open neighbourhood $\Nc^k$,
$\Nc^k \subset \Nc^{k-1}$,  such
that $\Pc(k)\cap\Nc^k$ is
dense in  $\Pc(k -1)\cap \Nc^k$. Since $\Pc(0)=\Om(\varepsilon,n+\delta)$,
Property (c) of Proposition \ref{period-mini} implies that
$\Pc(0)\cap\Nc^0=\Om^h_{n+\delta}\cap\Nc^0=\Nc^0$. These
   two properties show by recursion that $\Pc(k)\cap \Nc^k$ is dense in
   $\Nc^k$ for any $k$. Hence,
$\Om(n,n) \cap \Nc^{k_0}$ is dense in $\Nc^{k_0}$. As $f_0\in
\Om^{h}_{n}$ is arbitrary
and as $\Om^{h}_{n}$ is dense in $\Gm$, it follows that $\Om(n,n)$ is
dense in $\Gm$.

\begin{prop} \label{generic3A/2A}
For any $k \in\Nm$, $0 <k \leq k_0$, there exists a neighbourhood 
$\Nc^k \subset
\Nc^{k -1}$ of $f_0$ such that the set $\Pc(k)\cap \Nc^k$ is dense
in  $\Pc(k-1) \cap \Nc^k$.
\end{prop}

\begin{proof}
We apply Theorem \ref{A-th-Sard-Smale} from the
appendix with the following Banach spaces and functional $\Phi$.\HB
Let $0 <\delta^* \leq 1/4$, we set:
\begin{equation*}
\begin{split}
U = & \Big\{(T,u_0) \in
(\varepsilon,((3/2)^k+\delta^*)\varepsilon) \times H^s(S^1) \, |
\,
u_0 \notin {\overline {\Bc\Em}(f_0,n+\delta,r)}, 
\cr
&~~~\hbox{ and } \sup_{t \in
[0,((3/2)^k+\delta^*)\varepsilon]} 
\|S_{f_0}(t)u_0\|_{C^1(S^1)}
< n+\delta-\eta \sum_{j=1}^{k-1} j^{-4}- \frac{\eta}{2}k^{-4}\Big\}~, \cr
   Z  =  &H^s(S^1).
\end{split}
\end{equation*}
We would like to choose $\Pc(k-1)\cap \Nc^{k-1}$ as set $V$. However,
the space $\Gm$ is not even metrizable. To overcome this problem,
one has to work with functions defined on compact sets.
We notice that we can choose a neighbourhood $\Nc^k \subset
\Nc^{k-1}$ of $f_0$ such that, if
$$
\sup_{t \in
[0,((3/2)^k+\delta^*)\varepsilon]} \|S_{f_0}(t)u_0\|_{C^1(S^1)}
< n+\delta-\eta \sum_{j=1}^{k-1} j^{-4}- \frac{\eta}{2}k^{-4} ~,
$$
then, for any $g \in \Nc^{k}$,
$$
\sup_{t \in
[0,((3/2)^k+\delta^*)\varepsilon]} \|S_{g}(t)u_0\|_{C^1(S^1)}
< n+\delta-\eta \sum_{j=1}^{k-1} j^{-4}- \frac{\eta}{4}k^{-4}~.
$$
Let $R : g \in \Gm \mapsto Rg \in C^2(S^1 \times [-(n+2), n+2] \times
[-(n+2), n+2])$ be
the restriction operator defined by
$$
Rg =  g_{|S^1 \times [-(n+2), n+2] \times [-(n+2), n+2]}~.
$$
We set $V = R(\Pc(k-1)\cap \Nc^k)$ endowed with the topology of
$C^2(S^1 \times [-(n+2), n+2] \times [-(n+2), n+2],\Rm)$, which is a
separable Banach space.
The map $R$ is a continuous, open and surjective map.
Therefore the density of $R(\Pc(k)\cap \Nc^k)$ in $V$ is equivalent to
the density of  $\Pc(k)\cap \Nc^k$ in $\Pc(k-1)\cap \Nc^k$.\HB
We set $z=0$ and we consider the functional $\Phi : U \times V \to Z$
defined by
$$
\Phi : (T, u_0, f) \in U \times V \mapsto \Phi(T, u_0, f)= S_f(T)u_0 -
u_0~.
$$
As it will become clear below, $\Phi$ is a $C^2$-map from $U \times
V$ into $Z$.\\[2mm]

Let us make some remarks on the choice of this functional. First,
for any $(T,u_0)\in U$, the trajectory $S_g(t)u_0$ only depends on 
the value of $Rg$.
The use of the restriction operator $R$
does not affect any trajectory considered here.

  Next, we remark that $\Phi^{-1}(0) = \{ (T,u_0, f) \in U \times V \, | \,
S_f(t)u_0$ is a periodic orbit of \eqref{eq} of period $T \}$.
For any $(T, p(0), f) \in \Phi^{-1}(0)$, either $T$ is a simple
period of the periodic orbit $p(t)=S_f(t)p(0)$, or $T$ is not a
simple period. Since $f$ belongs to
$R(\Pc(k-1) \cap \Nc^k)$, any periodic orbit $p(t)$ of \eqref{eq} 
with period $T \in (0,(3/2)^{k -1}
\varepsilon]$ and such that $\sup_{t} \|p(t)\|_{C^1(S^1)} \leq n+\delta-\eta
\sum_{j=1}^{k-1} j^{-4}$ is hyperbolic. Therefore, as noticed in the 
remark following Proposition
\ref{simplehyper}, if $T$ is not a simple period of $p(t)$,
then $T$ must belong to $((3/2)^{k -1}\varepsilon, ((3/2)^{k 
}+\delta^*)\varepsilon]$ and
must be the minimal period of $p(t)$.

Working with sets involving a sequence of
  bounds of the type $n+\delta-\eta \sum_{j=1}^{k-1} j^{-4}-
\frac{\eta}{2}k^{-4}$ looks complicated and technical.
Actually, in the definition
of sets as $\Om^h_\mu$ or $\Om(A,\mu)$, we need large inequalities
(i.e. the symbols $\leq$), in order to obtain open sets, whereas
the definition of the open set $U$ requires strict inequalities (i.
e. symbols $<$). Thus, working at the same time with large and strict
inequalities, requires to introduce some intermediate bounds, which
look complicated and artificial.\\[2mm]

We now check that the hypotheses of Theorem \ref{A-th-Sard-Smale} are
  satisfied.\HB
Hypothesis i) of Theorem \ref{A-th-Sard-Smale} holds. Indeed, for any
$(T, p(0), f) \in \Phi^{-1}(0)$, we have
$$
D_{t,u}\Phi (T, p(0), f) (\tau, v) = p_t(T)\tau + (D_u (S_f(T)p(0))
-I)v \equiv p_t(T)\tau + (U_{f,p}(T,0) -I)v~,
$$
where the linearized operator $U_{f,p}$ has been introduced in
\eqref{linearized}. Since $U_{f,p}: H^s(S^1) \to H^s(S^1)$ is a
compact operator, the operator $U_{f,p}(T,0) -I$ is a Fredholm
operator of index $0$. We recall that
$\hbox{ dim (Ker }(U_{f,p}(T,0) -I)) \geq 1$.
Now two cases can occur. Either $p_t(T)$ belongs to
Im$(U_{f,p}(T,0) -I)$ and thus
$$
\hbox{ dim (Ker } D_{t,u}\Phi (T, p(0), f) )= \hbox{ dim (Ker
}(U_{f,p}(T,0) -I) )+1
$$
and
$$
\hbox{ codim (Im } D_{t,u}\Phi (T, p(0),
f) )= \hbox{ codim (Im}(U_{f,p}(T,0) -I)) .
$$
Or $p_t(T)$ does not belong to Im$(U_{f,p}(T,0) -I)$ and thus
$$
\hbox{ dim (Ker } D_{t,u}\Phi (T,p(0),f) )= \hbox{ dim (Ker
}(U_{f,p}(T,0) -I) )
$$ and
$$
\hbox{ codim (Im } D_{t,u}\Phi (T, p(0),f) )= \hbox{ codim
(Im}(U_{f,p}(T,0) -I)) -1 .
$$
In both cases, we
conclude that $D_{t,u}\Phi (T, p(0), f)$ is a Fredholm operator of
index $1$. We notice that, if $D_{t,u}\Phi (T, p(0), f)$ is a
surjective map from $\Rm \times H^s(S^1)$ into itself, then the
dimension of $\hbox{ Ker } D_{t,u}\Phi (T, p(0), f))$ is equal to
$1$.  We also remark that, if $T$ is a simple period of $p(t)$, then
$D_{t,u}\Phi (T, p(0), f)$ is a surjective map.  In particular, for
any $(T, p(0), f) \in \Phi^{-1}(0)$ such that the minimal period of
$S_f(t)p(0)$ is less than (or equal to) $(3/2)^{k -1}\varepsilon$,
we know that $D_{t,u}\Phi (T, p(0), f)$ is a surjective map, since $f$
belongs to $\Pc(k-1)$.

\vskip 2mm

\noindent We next show that Hypothesis ii) of Theorem
\ref{A-th-Sard-Smale} is also satisfied.  By the above considerations,
we are reduced to proving that $D_{t,u,f}\Phi (T, p(0), f)$ is a
surjective map, only in the case where $T \in [(3/2)^{k
-1}\varepsilon, ((3/2)^{k }+\delta^*)\varepsilon]$ is the minimal
period of $p(t)=S_f(t)p(0)$.  An easy computation shows that
$$
D_{t,u,f}\Phi (T, p(0), f) (\tau, v,g) = p_t(T)\tau + (U_{f,p}(T,0)
-I)v + \Sigma_{f,p}(T)g~,
$$
where $\Sigma_{f,p}(t)g =w(t)$ is the solution of the following affine
equation
\begin{equation}
\label{eqg}
\begin{split}
w_t(x,t) &=w_{xx}(x,t)+ D_{u}f(x,p,p_{x})w
+ D_{u_x}f(x,p,p_{x})w_x + g(x,p,p_{x})~,  \cr
w(x,0) & =0~.
\end{split}
\end{equation}
We remark that
\begin{equation}
\label{eqgBIS}
w(t) = \int_{0}^{t} U_{f,p}(t,s)g(\cdot, p(s),p_x(s)) ds~.
\end{equation}
The map $D_{t,u,f}\Phi (T, p(0), f)$ is surjective if and only if,
for any $h$ in
$H^s(S^1)$, there exists $(\tau, v,g) \in \Rm \times H^s(S^1) \times
R\Gm$, such that,
$$
p_t(T)\tau + (U_{f,p}(T,0) -I)v + \Sigma_{f,p}(T)g = h~.
$$
Due to the Fredholm alternative, $h - \Sigma_{f,p}(T)g$ belongs to
the image of $U_{f,p}(T,0) -I$ if and only if
\begin{equation}
\label{integrale0}
\int_{S^1} (h - \Sigma_{f,p}(T)g) \varphi^*(x) dx = 0~,
\end{equation}
for any solution $\varphi^*$ of the adjoint equation
\begin{equation}
\label{adjoint}
(U_{f,p}(T,0))^*\varphi^*  =\varphi^*~.
\end{equation}
Let $\varphi^*_i$, be a (at most two-dimensional) basis of $\hbox{
Ker }(U_{f,p}(T,0)-Id)^*$. We must find $g$ such that
$\int_{S^1} (\Sigma_{f,p}(T)g) \varphi^*_i (x) dx = \int_{S^1}  h(x)
\varphi^*_i (x) dx$, for any $i$. The surjectivity of the map
$g \mapsto (\int_{S^1} (\Sigma_{f,p}(T)g) \varphi^*_i (x) dx )_i$ is
equivalent to the non-existence of $(c_i)_i$ such that $\sum_{i}c_i
\int_{S^1} (\Sigma_{f,p}(T)g) \varphi^*_i (x)dx =0$, for every
$g  \in R\Gm$. Thus, we are
reduced to proving that there is no solution $\varphi^* \ne
0$ of the adjoint equation $(U_{f,p}(T,0))^*\varphi^* =\varphi^*$
such that
$\int_{S^1} (\Sigma_{f,p}(T)g) \varphi^*(x) dx$ $ = 0$, for every
$g  \in R\Gm$. In other terms, we are reduced to prove that, for any
solution $\varphi^* \ne 0$ of the adjoint equation
$(U_{f,p}(T,0))^*\varphi^* =\varphi^*$, there exists $g \in R\Gm$ such
that
\begin{equation}
\label{nonnulle}
\int_{S^1} (\Sigma_{f,p}(T)g) \varphi^*(x) dx \ne 0~.
\end{equation}
Using the property \eqref{eqgBIS}, we see that the condition
\eqref{nonnulle} is equivalent to the fact that, for any $\varphi^*$,
there exists $g \in R\Gm$ such that
\begin{equation}
\label{surjection1}
\begin{split}
\int_{S^1}\int_{0}^{T} U_{f,p}(T,s)(g(x, &p(x,s),p_x(x,s)))
\varphi^*(x) ds dx \cr
= &\int_{S^1}\int_{0}^{T} g(x, p(x,s),p_x(x,s))
[(U_{f,p}(T,s))^*\varphi^*](x) ds dx  \ne 0~.
\end{split}
\end{equation}
Since $\varphi^*\neq 0$, there exist $x_0 \in S^1$ and $t_0 \in
[0,T)$ such that $[(U_{f,p}(T,t_0))^*\varphi^*](x_0) \ne 0$. It is easy
to construct a
regular bump function $ g(x,u,u_x)$ which vanishes outside a small
neighborhood of $(x_0,p(x_0,t_0), p_x(x_0,t_0))$ and is positive in
this neighborhood.
Due to the injectivity property of Lemma \ref{injectif}, for $x_0$
fixed, there exists no other time $t_1 \in [0,T)$ such that
$p(x_0,t_1)=p(x_0,t_0)$ and
$p_x(x_0,t_1)=p_x(x_0,t_0)$. Therefore, the function
$(x,s)\longmapsto g(x, p(x,s),p_x(x,s))$ is a regular bump function
concentrated around $(x_0,t_0)$. For such a choice of $g$, the
condition \eqref{surjection1} is thus satisfied.\\[2mm]

Since all the hypotheses of Theorem \ref{A-th-Sard-Smale} hold, there
exists a generic subset $V_1$ of $V$ such that, for any $f \in V_1$,
the map $D_{t,u}\Phi (T, p(0), f)$ is a  surjective map. If $p_t(T)$
belongs to
Im$(U_{f,p}(T,0) -I)$, then $U_{f,p}(T,0) -I$ is surjective, which is
not possible since it is a Fredholm operator of index $0$ and that its
      kernel contains $p_t(0)$. Thus, $p_t(T)$ does not belong to
Im$(U_{f,p}(T,0) -I)$ and $1$ is an eigenvalue of algebraic
      multiplicity $1$. Then, $\hbox{ dim (Ker (}U_{f,p}(T,0) -I))=
\hbox{ codim (Im (}U_{f,p}(T,0) -I))=1$ and thus $1$ is an
eigenvalue of geometric multiplicity $1$. The proposition is then
proved. \end{proof}


\section{Generic hyperbolicity for nonlinearities
independent of $x$}\label{Sec-inde}
The purpose of this section is the proof of Theorem \ref{th-inde-x}.
We denote by $\Gm^i$ the subspace of $\Gm$ consisting of functions
independent of $x$. We use the terminology of Proposition
\ref{th-Ang-Fied}. Notice that if $u$ is a frozen wave, then every
spatial translation
$u(\cdot -x_0)$ of $u$ is also a frozen wave.  This means that a frozen
wave always belongs to a circle of equilibria and is never a
hyperbolic equilibrium of \eqref{eq-inde-x}.  However, it is more
natural to consider the frozen waves as particular cases of rotating
waves with speed $c=0$.  For the remaining part of this section, we
group the waves with
speed $c\in\Rm$ into a single category and simply call them ``waves"
if there is no need to distinguish between rotating and frozen
waves.  To simplify the notations, we
say that a frozen wave $u$ is hyperbolic if the corresponding circle
of equilibria is normally hyperbolic, that is, if $0$ is a simple
eigenvalue of the linearized operator with the eigenfunction $u_x$ and
is the only eigenvalue with real part equal to $0$.\\

     If $u(x,t)=v(x-ct)$ is a wave solution of \eqref{eq-inde-x} of speed
$c\in\Rm$,
     then $v$ is an equilibrium of the equation
     \begin{equation}
     \label{eqv}
w_t= w_{xx} +  f(w,w_x) + c w_x~.
\end{equation}
For this reason, the first step of the proof of Theorem
\ref{th-inde-x} consists in eliminating the time dependence in the
problem. Let
$L:H^{2}(S^1)\longrightarrow L^2(S^1)$ to be the linearized operator
defined by
$$ L\varphi=\varphi_{xx}+f'_u(v,v_x)\varphi+f'_{u_x}(v,v_x)
\varphi_x+c\varphi_x~.$$
\begin{lemma}\label{lemme-inde-1}
A wave $u(x,t)=v(x-ct)$ of \eqref{eq-inde-x} is hyperbolic if and only
if $0$ is a simple eigenvalue of $L$. Moreover, if $v$ is spatially
$\frac1n$-periodic, i.e. $t\mapsto u(t)$ is $\frac1{nc}$ periodic,
then any solution $\psi$ of $L\psi=0$ or of $L^*\psi=0$ is also
spatially periodic of period $\frac 1n$.
\end{lemma}

\begin{proof}
If $u$ is a frozen wave, then the first assertion of Lemma
\ref{lemme-inde-1} is a direct consequence of Proposition
\ref{spectreLe}.\\
Let $u(x,t)=v(x-ct)$ be a rotating wave of period $T=1/c$.  We
consider the operator $U(T,0):H^s(S^1)\longrightarrow H^s(S^1)$
defined by $U(T,0)\varphi_0=\varphi(T)$ where $\varphi(t)$ is the
solution of
\begin{equation}\label{eq-lemme-inde-1}
\left\{\begin{array}{l}\varphi_t=\varphi_{xx}+f'_u(u,u_x)\varphi+f'_{u_x}(u,u_x)\varphi_x~,

\\
\varphi(x,0)=\varphi_0(x)~.\end{array}\right.
\end{equation}
We set $\varphi(x,t)=\psi(x-ct,t)$, so that $\psi$ satisfies $\psi_t=L
\psi$. Let $V(T,0)=e^{TL}$. Since $cT=1$, there is a perfect
correspondence between the spectrum of $U(T,0)$ and the spectrum of
$V(T,0)$. In particular, Proposition \ref{simplehyper} shows that the
rotating wave $u$ is hyperbolic if and only if $1$ is a simple
eigenvalue of $V(T,0)$. Notice that $v_x \ne 0$ belongs to
$\hbox{Ker }L$. Thus, by Proposition \ref{spectreLe}, zero is the only
eigenvalue of zero real part for $L$. The spectral theorem recalled
in \eqref{spectralmap} implies that $u$ is hyperbolic if and only if
$0$ is a simple eigenvalue of $L$.\\
Assume now that $v$ is $\frac1n$-periodic. Let $\tilde L$ be the
restriction of $L$ to $H^s_{1/n-per}(S^1)=H^s(\Rm/(1/n)\Zr)$ and let
$\tilde v_x$ the
restriction of $v_x$ to $\Rm/(1/n)\Zr$.  We have $\tilde L\tilde
v_x=0$, which implies that $0$ is also an eigenvalue of $\tilde L^*$.
By Proposition \ref{spectreLe},
the only eigenvalue of $\tilde L$ or $\tilde L^*$ with zero real part
is zero. Using the remark following Proposition
\ref{simplehyper}, we prove that every solution of $L\psi=0$ (resp.
$L^*\psi=0$) is necessarily $\frac 1n-$periodic and solution of $\tilde
L\tilde \psi=0$ (resp. $\tilde L^*\tilde \psi=0$). \end{proof}

\noindent We have to consider the case of homogeneous equilibria
separately.
\begin{lemma}\label{lemme-inde-3}
There exists a generic dense subset $\Om^{hom}$ of $\Gm^i$ such that
every homogeneous equilibrium point of \eqref{eq-inde-x} is
hyperbolic.
\end{lemma}
\begin{proof}
The lemma is a direct consequence of the following computation.  If
$u\in\Rm$ is a homogeneous equilibrium, then there exist a function
$\varphi\in H^2(S^1)$ and $\mu\in\Rm$ such that
$\varphi_{xx}+f'_u(u,0)\varphi+f'_{u_x}(u,0)\varphi_x=i\mu\varphi$ if
and only if there exists $k\in\Zr$ such that $2kf'_{u_x}(u,0)\pi=\mu$
and $f'_u(u,0)=4 k^2\pi^2$.  Thus
$\Om^{hom}=\{f\in\Gm^i \, | \, f(u,0)=0\Rightarrow f'_u(u,0)\not\in 4\pi^2\Nm
\}$ is a suitable choice.
\end{proof}

The main step of the proof of Theorem \ref{th-inde-x} consists in
applying Sard-Smale Theorem to deal with the waves.
\begin{lemma}\label{lemme-inde-4}
There exists a generic subset $\Om^{wav}$ of $\Gm^i$ such that every
wave of \eqref{eq-inde-x} is hyperbolic.
\end{lemma}
\begin{proof}
The proof is very similar to the one of Proposition
\ref{generic3A/2A}.  The main tool is Theorem \ref{A-th-Sard-Smale} of
the appendix.  We recall that $\Gm^i$ is not a metrizable space.  To
overcome this difficulty, we write
$\Om^{wav}=\cap_{n\in\Nm}\Om^{wav}_n$ where $\Om^{wav}_n$ is the set
of functions of $\Gm^i$ such that every wave $v$ of \eqref{eq-inde-x}
with $\|v\|_{C^1}<n$ is hyperbolic.  We next prove that, for each
$n\in\Nm$, $\Om^{wav}_n$ is a generic subset of $\Gm^i$, which implies
that $\Om^{wav}$ is a generic subset of $\Gm^i$, since it is a
countable intersection of generic sets.\HB
Let $n\in\Nm$.
We again introduce the restriction operator $R: g \in \Gm^i \mapsto
Rg \in C^2([-(n +2), n+2] \times [-(n +2), n+2],\Rm)$ defined by
$$ Rg =  g_{|[-(n+2), n +2] \times [-(n+2), n+2}]~.$$
We
remind
that $R\Gm^i=C^2([-(n +2), n+2] \times [-(n +2), n+2],\Rm)$ is a
separable Banach space.
As in the proof of Proposition \ref{generic3A/2A}, it is sufficient to work
with $C^2([-(n +2), n+2] \times [-(n +2), n+2],\Rm)$ since the values 
of $f$ outside
$[-n, n] \times [-n, n]$ do not matter for the genericity of the set 
$\Om^{wav}_n$.

We apply Theorem \ref{A-th-Sard-Smale} of the appendix with the following
spaces and functional.  We set $U=\{v\in H^{2}(S^1) \, | \, v_x \not \equiv
0, \|v\|_{C^1}<n\}\times\Rm$, $V=R\Gm^i$ and $Z=L^2(S^1)$.  We set $z=0$ and
$$\Phi:~\left(\begin{array}{ccc} U\times V & \longrightarrow & Z \\
(v,c,f)&\longmapsto & v_{xx}+f(v,v_x)+cv_x \end{array}\right)~.$$
First notice that Hypothesis iii) of Theorem \ref{A-th-Sard-Smale} is
satisfied and that a point $(v,c,f)$ belongs to $\Phi^{-1}(0)$ if and
only if $u(x,t)=v(x-ct)$ is a wave of speed $c$ for \eqref{eq-inde-x}.
We set
$$
L\varphi=\varphi_{xx}+f'_u(v,v_x)\varphi+f'_{u_x}(v,v_x)\varphi_x+c\varphi_x~.
$$
We recall that, as proved in Lemma \ref{lemme-Le-Fredholm}, $L$ is a
Fredholm operator of index 0.\\
Hypothesis i) of Theorem \ref{A-th-Sard-Smale} holds.  Indeed,
$D_{v,c}\Phi(v,c,f).(\varphi,d)=L\varphi+d v_x$. As in the proof of
Proposition \ref{generic3A/2A}, there are two cases. Either $v_x$
belongs to Im$(L)$ and thus
dim(Ker$(D_{v,c}\Phi)$)=dim(Ker$(L)$)+1,  codim$ $(Im$ (D_{v,c}\Phi)
$)=codim(Im$(L)$).  Or $v_x$ is not in Im$(L)$ and then dim
(Ker$ (D_{v,c}\Phi) $)=dim(Ker$(L)$) and
codim(Im$(D_{v,c}\Phi)$)= codim(Im$(L)$)-1.  In both cases, since $L$
is Fredholm of index 0, $D_{v,c}\Phi$ is a Fredholm
operator of index $1$.\\
Hypothesis ii) of Theorem \ref{A-th-Sard-Smale} is also satisfied.  Indeed,
let $h\in L^2(S^1)$.  We have to find $v$ and $g$ such that
$D\Phi(v,c,f).(\varphi,0,g)=L\varphi+g(v,v_x)=h$.  Due to the Fredholm
alternative, $h-g(v,v_x)$ belongs to Im($L$) if and only if
$\int_{S^1}(h-g(v,v_x))\psi=0$ for all $\psi$ solution of $L^*\psi=0$.
Let $(\psi_i)$ be a (at most two-dimensional) basis of Ker($L^*$), we
must find $g$ such that $\int g(v,v_x)\psi_i=\int h\psi_i$ for all
$i$.  The surjectivity of the map $g\mapsto \left(\int
g(v,v_x)\psi_i\right)_i$ is equivalent to the non-existence of
$(c_i)\neq 0$ such that $\sum c_i\int g(v,v_x)\psi_i=0$ for all $g$.
Thus, setting $\psi=\sum c_i\psi_i$, we are reduced to prove that
there is no $\psi\neq 0$ such that $L^*\psi=0$ and
$\int_{S^1}g(v,v_x)\psi=0$ for all $g\in R\Gm^i$.
Let $\psi$ satisfying
$L^*\psi=0$ and $\int_{S^1}g(v,v_x)\psi=0$ for
all $g\in R\Gm^i$.  Let
$m \ne 0$ be the integer such that $1/m$ is the minimal period of the
function $v$.  By Lemma \ref{lemme-inde-1}, $\psi$ is also
$1/m-$periodic and $\int_{S^1}g(v,v_x)\psi=0$ is equivalent to
$\int_0^{1/m}g(v,v_x)\psi=0$.  Since $v$ is a periodic solution of
minimal period $1/m$ of a second
order ordinary differential equation, the map $x\in [0,1/m)\mapsto
(v(x),v_x(x))\in\Rm^2$ is injective and, hence
$\int_0^{1/m}g(v,v_x)\psi=0$ for all $g$ implies that $\psi\equiv0$.
This shows the surjectivity of $D\Phi(v,c,f)$.\\
All the assumptions of Theorem \ref{A-th-Sard-Smale} being satisfied,
there exists a generic subset $R\Om^{wav}_n$ of $R\Gm^i$ such that for any
$f\in R\Om^{wav}_n$ and for any wave $u(x,t)=v(x-ct)$ of
\eqref{eq-inde-x} in $U$, the map $D_{v,c}\Phi:~(\varphi,d)\mapsto
L\varphi+dv_x$ is surjective.  If $v_x$ belongs to Im($L$), this
surjectivity implies that $L$ is surjective, which is not possible
since $L$ is a Fredholm operator of index 0 and since Ker($L$) is not
$\{0\}$ because it contains $v_x\neq 0$.  Thus, $v_x$ does not belongs
to Im($L$), which means that $0$ is an algebraically simple
eigenvalue for $L$.  Moreover, dim(Ker$L$)=codim(Im($L$))=1 which
means that $0$ is a geometrically simple eigenvalue of $L$.  The
conclusion is then given by Lemma \ref{lemme-inde-1}.\end{proof}

{\noindent \bf Proof of Theorem \ref{th-inde-x}: }Applying Lemmas
\ref{lemme-inde-3}
and \ref{lemme-inde-4}, we obtain a generic subset of $\Gm^i$ such
that all homogeneous
equilibria and all waves of \eqref{eq-inde-x} are hyperbolic. To
obtain Theorem \ref{th-inde-x},
it only remains to remove the frozen waves. Indeed, we recall
that a frozen wave is never
hyperbolic as a non-homogeneous
equilibrium point of
\eqref{eq-inde-x} (and not as a wave as it was considered in the
above lemmas).
As already noticed in \cite{FRW}, if we replace $f(u,u_x)$ by
$f(u,u_x)+\varepsilon u_x$, we can ``unfreeze'' every frozen wave. To
be sure that we are not ``freezing'' some rotating wave, we remark
that there is at most a countable number of hyperbolic waves. Indeed,
by arguments similar to the ones used in the proof of Property (d) of
Proposition \ref{period-mini}, one can show that there is a finite
number of periodic orbits $p$ of period less that $n$ and satisfying
$\sup_t \|p(t)\|_{\Cm^1}\leq n$. Thus, we can choose $\varepsilon>0$
as small as
wanted to unfreeze the frozen waves without freezing any rotating wave.
\HB We emphasize that, for a general
system, without the constraint
on the period, the number of hyperbolic periodic orbits could be
infinite, even in bounded sets, since hyperbolic
periodic orbits might pile up on a homoclinic orbit of a hyperbolic
equilibrium point. In the case of Equation \eqref{eq-inde-x}, one
should be able to remove the constraint on the period, since there do
not exist homoclinic orbits (see \cite{FRW}).
{\hfill $\square$\\}


\renewcommand{\thesection}{A}
\section{Appendix : Sard-Smale theorem}

Let $M$ and $M'$ be two differentiable Banach manifolds and let
$f:M\longrightarrow M'$ be a differentiable map.  We say that $y\in
M'$ is a regular value of $f$ if, for any $x\in M$ such that $f(x)=y$,
the differential $Df(x):T_xM\longrightarrow T_yM'$ is surjective.  The
points of $M'$ which are not regular are said to be critical.  The classical
theorem of Sard says that, if $U$ is an open set of $\Rm^p$ and if
$f:U\longrightarrow \Rm^q$ is of class $\Cm^s$ with $s>max(p-q,0)$,
then, the set of critical values of $f$ in $\Rm^q$ is of Lebesgue
measure zero. Using Fredholm operators and a Lyapunov-Schmidt method,
Smale has generalized the Sard Theorem to infinite-dimensional spaces
(for Fredholm operators, we refer to \cite{Bonic} for
example). The Sard-Smale Theorem stated below is an application of
the Smale Theorem.

\begin{theorem}\label{A-th-Sard-Smale}
Let $X,Y,Z$ be three smooth Banach manifolds. Let $U\subset X$, $V\subset Y$
be two open
sets, $\Phi:U\times V\longrightarrow Z$ be a map of class $\Cm^k$
($k\geq 1$) and $z$ be a
point of $Z$. \\
We assume that:
\begin{enumerate}
    \item[i)] $\forall (x,y)\in \Phi^{-1}(z), D_x\Phi(x,y)$ is a Fredholm
       operator of index strictly less than $k$,
    \item[ii)] $\forall (x,y)\in \Phi^{-1}(z),D\Phi(x,y)$ is surjective,
    \item[iii)]  $X$ and $Y$ are separable.
\end{enumerate}
Then $\Theta=\{y\in V \, | \, z$ is a regular value of $\Phi(.,y)\}$ is a
generic subset of $V$.
\end{theorem}
The proof of Theorem \ref{A-th-Sard-Smale} can be found in
\cite{sauttemam} or \cite{Quinn}
(for stronger versions, see also \cite{He}). We also refer to
\cite{RJ} where one can find a short review on PDE
generic results, as well as
adaptations of the transversality theorems to the notion of
prevalence, which is a notion
of ``almost always'' different from the genericity.\\

\bibliographystyle{amsplain}

\end{document}